\documentclass{article}
\usepackage{amsfonts, amsbsy, amssymb}

\newtheorem{lemma}{Lemma}

\newtheorem{theorem}{Theorem}
\def\F{{\mathbb F}}

\title{An improved estimate on sums of product sets\footnotetext{AMS subject classification 11T, 52C} }

\author{Misha Rudnev\thanks{University of Bristol, Bristol BS8 1TW UK, {\sf
m.rudnev@bris.ac.uk}}}
\begin{document}
\maketitle

\begin{abstract} In a recent paper \cite{Gl} A. Glibichuk proved that if $A,B$ are subsets of an arbitrary finite filed $\F_q$, such that $|A||B|>q$, then $16AB = \F_q$. We improve this to
$10AB = \F_q.$
\end{abstract}

\renewcommand{\theenumi}{\roman{enumi}}

Let $\F_q$ be a finite field of  $q$ elements. As usual, for $A,B\subseteq \F_q,\,\xi\in \F_q^*$, where $\F_q^*$ stands for the multiplicative group of $\F_q$,  we denote
$$\begin{array}{cccc}
A+B=\{a+b:\,a\in A,\,b\in B\}, & & AB=\{ab:\,a\in A,\,b\in B\},\\ \hfill \\ \xi B=\{\xi\}B,\;\;\;\;-A=\{-1\}A,&& dA = \underbrace{A+\ldots+ A}_{d \mbox{ {\tiny times}}},\end{array}$$
for $d\in \mathbb N$. The main results of this note are as follows.
\begin{theorem} Let $A,B\subseteq\F_q$, such that the product of their cardinalities
$|A||B|>q$. Then $10AB=\F_q$. \label{th}\end{theorem}
Theorem \ref{th} is an improvement on a recent result of Glibichuk  (\cite{Gl}) who showed $16AB=\F_q$ as a consequence of a stronger claim $8AB=\F_q$ if one of the sets $A,B$ is symmetric or antisymmetric (which also implies that $8AB=\F_q$ as long as $|A||B|\geq 2q$.) Our claim is based on one simple observation and a slightly more elaborate use of symmetry. The constant $10$ regarding $10AB=\F_q$ is unlikely to be optimal. A more general question is: under the assumption $|A||B|>(1+c)q$, what is the smallest integer $d(c)$ so that $dAB$ cover a fraction $\frac{1}{C(c)}$ of elements of $\F_q$? One is tempted to believe that $d=2$ should not generally suffice for $c=o(1),\, C=O(1)$ yet we are unaware of constructive evidence to this. D. Hart and A. Iosevich (\cite{HI}) conjectured that in the case $A=B$, the condition $|A|\geq C_\epsilon q^{\frac{1}{2}+\epsilon}$ should suffice for $2A^2$ to cover the whole of $\F_q$.

Geometrically the integer $d$  has the meaning of the dimension, so that $d$-dimensional Cartesian products $A_d$, $B_d$ of $A$ and $B$, respectively, with itself generate sufficiently many distinct dot products $\boldsymbol a\cdot \boldsymbol b=a_1b_1+\ldots+a_db_d$, where $\boldsymbol a = (a_1,\ldots,a_d) \in A_d$, $\boldsymbol b = (b_1,\ldots,b_d)\in B_d$. This geometric interpretation of the arithmetic problem was used quite elegantly in \cite{HI} by way of Fourier analysis; it also puts this problem under the shibboleth of ``hard Erd\"os problems'' and arguably distinguishes the sum $AB+AB+\ldots$, within the more general case when some of the plus signs are replaces by minuses. It follows in particular from Theorem  \ref{th} that in $d=10$ the set of dot products of elements of $A_d$ with itself is the whole field $\F_q$, as long as $|A_d|>q^{\frac{d}{2}}$, quite in the spirit of, say, the Erd\"os-Falconer distance problem. (In the case $A=B$, the following proof can be interpreted that there exists $\boldsymbol a_*\in A_d$, such that dot products with $\boldsymbol a_*$ cover $\F_q$.)

\medskip
The proof of Theorem \ref{th} follows from massaging the results of Lemmas 1,2, and 4 in \cite{Gl}, which we formulate in a form suitable for immediate use in the sequel.
The following lemma the key point of the argument. It has recently appeared in arithmetic combinatorics literature on numerous occasions, following up on a ``statement about generic projections'' by Bourgain, Katz, and Tao (\cite{BKT}, Lemma 2.1).

\medskip
By default, the sets $A,B$ in the sequel always satisfy the conditions of Theorem \ref{th}.

\addtocounter{lemma}{1}
\begin{lemma}  [\cite{Gl}, {\small Lemmas 1,2}] \label{Glib} There exists $\xi\in \F_q^*$, such that \begin{enumerate} \item the equation
\begin{equation}
a_1+\xi b_1 = a_2+\xi b_2
\label{deq}\end{equation}
has strictly fewer than ${\displaystyle \frac{2}{q} |A|^2|B|^2}$ solutions $(a_1,b_1,a_2,b_2)\in A\times B\times A\times B,$ \item both sets
$C^\pm_\xi = A\pm\xi B$ are such that $$|C^\pm_\xi|>\frac{q}{2}.$$\end{enumerate} \end{lemma}

The statement (ii) follows form (i) by Cauchy-Schwartz inequality. The reason why there are two sets $C^\pm_\xi$ for the same $\xi$ is that the equation (\ref{deq}) can be rewritten as
$$a_1-\xi b_2 = a_2-\xi b_1.$$

As a separate term, used in the sequel, let us call $y$  in $C^+_\xi$  ($C^-_\xi$) {\em involved} if it allows for more than one representation $y=a+\xi b$ ($y=a-\xi b$) in terms of elements of $(A,B)$. Since $|A||B|>q$, there exists an involved $y$ in $C^+_\xi$  ($C^-_\xi$).

\medskip
The second pre-requisite we need is as follows.

\begin{lemma} [\cite{Gl}, {\small Lemma 4}] \label{two} If $C\subset\F_q$ is such that $|C|>\frac{q}{2}$, then $2C=\F_q$.
\end{lemma}

We now turn to the proof of Theorem \ref{th}. Let $\xi,\,C^\pm_\xi$ come from Lemma \ref{Glib} and be fixed once and for all.
\begin{lemma} \label{pr1} If there exists $a\in -A$ but not in $A+A$ or there exists $b\in -B$, but not in $B+B$, then $10AB =\F_q$.\end{lemma}
{\sf Proof:}  Fix $a\in -A$ which is not in $A+A$ (if there is no such $a$, then there is $b\in -B$ but not in $B+B$ and we swap $A$ and $B$). By Lemmas \ref{Glib} (ii) and \ref{two}, we have $2C^+_\xi = \F_q$. Therefore we have
$$
a= (a_1+a_2) + \xi(b_1+b_2),
$$
for some $a_1,a_2\in A$ and  $b_1,b_2\in B$, further fixed. Since $a\neq a_1+a_2$, this unambiguously determines
$$
\xi = - \frac{a_1+a_2-a}{b_1+b_2} = - \frac{a_1+a_2+a_3}{b_1+b_2},
$$
where $a_3=-a\in A$. Hence, by Lemma \ref{Glib} (ii), the set
$$(b_1+b_2)C^-_\xi = \{(b_1+b_2)a+ (a_1+a_2+a_3)b:\,(a,b)\in A\times B\}$$ has cardinality in excess of $\frac{q}{2}$ and is clearly a subset of $5AB$.  Lemma \ref{pr1} now follows by Lemma \ref{two}. $\Box$

\medskip
In view of Lemma \ref{pr1} we may now assume that all elements of $-A$ and $-B$ belong to $A+A$ and $B+B$, respectively. Each pair of the four sets $C^\pm_\xi,\,-C^\pm_\xi$ intersects, because each set has cardinality greater than $\frac{q}{2}$. Let $x$ belong to the intersection of some different two of those four sets; in the sequel a pair of those sets means a pair of distinct ones. Let us call $x$ {\em trivial} if it does not enable one to determine $\xi$ unambiguously, otherwise it is non-trivial. For instance, if $x\in C^+_\xi\cap - C^-_\xi$ we have, for some $(a_1,b_1,a_2,b_2)\in A\times B\times A\times B$,
$$
a_1+\xi b_1 = -a_2+\xi b_2,
$$ and $x$ will be trivial if all its such representations have $a_1=-a_2$.

\begin{lemma} Theorem \ref{th} follows if there exists a non-trivial $x$ in some pair intersection of the sets $C^\pm_\xi,-C^\pm_\xi$. \label{pr2}\end{lemma}
{\sf Proof:} Suppose, there is a non-trivial $x\in C^+_\xi\cap -C^+_\xi$. Then
$$
x= a_1+\xi b_1 = -a_2 -\xi b_2,
$$
for some $(a_1,b_1,a_2,b_2)\in A\times B\times A\times B,$ and $b\neq -b_2$. So
$$
\xi = \frac{a_1+a_2}{b_1+b_2},
$$
with non-zero denominator. Hence, by Lemma \ref{Glib} (ii), the set
$$(b_1+b_2)C^+_\xi = \{(b_1+b_2)a+ (a_1+a_2)b:\,(a,b)\in A\times B\}$$ has cardinality in excess of $\frac{q}{2}$ and is clearly a subset of $4AB$; then Lemma \ref{two} ensures that
 $8AB$ covers $\F_q$.

Suppose now a non-trivial $x$ lives in $C^+_\xi\cap C^-_\xi$. Then
$$
x= a_1+\xi b_1 = a_2 -\xi b_2,
$$
for some $(a_1,b_1,a_2,b_2)\in A\times B\times A\times B,$ and $b\neq -b_2$. So
$$
\xi = \frac{a_2-a_1}{b_1+b_2} = \frac{a_2+a_3+a_4}{b_1+b_2},
$$
for some $a_3, a_4\in A$ (by the assumption $-A\subseteq A+A$) and with non-zero denominator. Hence, by Lemma \ref{Glib} (ii), the set
$$(b_1+b_2)C^+_\xi = \{(b_1+b_2)a+ (a_1+a_2+a_3)b:\,(a,b)\in A\times B\}$$ has cardinality in excess of $\frac{q}{2}$ and is clearly a subset of $5AB$; by Lemma \ref{two} now $10AB$ covers $\F_q$.

Similarly, if a non-trivial $x$ lives in $C^+_\xi\cap -C^-_\xi$, we have
$$
x= a_1+\xi b_1 = -a_2 +\xi b_2
$$
for some $(a_1,b_1,a_2,b_2)\in A\times B\times A\times B,$ and $b_1 \neq b_2$. So
$$
\xi = \frac{a_1+a_2}{b_2-b_1} = \frac{a_1+a_2}{b_2+b_3+b_4},
$$
for some $b_3, b_4\in B$ (by the assumption $-B\subseteq B+B$) and with non-zero denominator. Hence, by Lemma \ref{Glib} (ii), the set
$$(b_2+b_3+b_4)C^+_\xi = \{(b_2+b_3+b_4)a+ (a_1+a_2)b:\,(a,b)\in A\times B\}$$ has cardinality in excess of $\frac{q}{2}$ and is clearly a subset of $5AB$; by Lemma \ref{two} now $10AB$ covers $\F_q$.

Other three pairs out of the four sets $C^\pm_\xi,\,-C^\pm_\xi$ are clearly amenable to one of the three cases above. $\Box$

\medskip
It remains to establish that a non-trivial $x$ exists. Let $\bar A$ denote a symmetric part of $A$ and $\tilde A$ its antisymmetric part. I.e.
$$
\bar A =\{a\in A:\,-a\in A\},\qquad \tilde A = A\setminus \bar A,
$$
the same for $B$. Suppose, a non-trivial $x$ does not exist. This entails two consequences. Firstly, all the sets
\begin{equation}\label{sets}
\begin{array}{llllllll} \tilde A+\xi \tilde B,&\tilde A-\xi \tilde B,& -\tilde A+\xi \tilde B,&-\tilde A-\xi \tilde B,\\ \hfill \\
\bar A+\xi \tilde B,&\bar A-\xi \tilde B,& \tilde A+\xi \bar B,&-\tilde A+\xi \bar B,\\ \hfill \\
\bar A+\xi\bar B
\end{array}\end{equation}
must be pairwise-disjoint, or there would exist a non-trivial $x$. Secondly, the diophantine equation (\ref{deq}) restricted to {\em any one of} the last five of these sets may not have any {\em involved} solutions. Indeed, if one of these two conditions failed, then using  $\bar A\cap\tilde A=\emptyset$, as well as $\tilde A\cap\,-\tilde A=\emptyset$, and $\bar A=-\bar A$, the same for $B$, it would be possible to express $\xi$ unambiguously in the form either ${\displaystyle \frac{a_1+a_2}{b_1+b_2}}$, or ${\displaystyle \frac{a_1-a_2}{b_1+b_2}}$, or ${\displaystyle \frac{a_1+a_2}{b_1-b_2}}$ and then repeat the reasoning inside the proof of Lemma \ref{pr2}. The only doubling, i.e the existence of involved solutions of (\ref{deq}), may occur when $a_1,a_2,b_1,b_2$ appearing therein are restricted to one of the first four sets on the disjoint list (\ref{sets}). The fact that $x$ is involved would then yield an unambiguous expression ${\displaystyle \xi = \frac{a_1-a_2}{b_1-b_2}}$, but the latter defies the argument inside the proof of lemma \ref{pr2} (i.e. the same argument would only result in $12AB=\F_q$).

It remains to bring the reasoning to absurd. Let $|\tilde A|=u|A|,$ $|\tilde B|=v|B|$ for some $u,v\in (0,1)$. (If one of the sets is symmetric or antisymmetric, i.e. $u$ or $v$ is $0$ or $1$, then by Glibichuk's result (\cite{Gl}) $8AB$ covers $\F_q$.)

Let us estimate from below the size of the first four sets on the list (\ref{sets}). Lemma \ref{Glib} (i) provides the upper estimate for the total number of solutions of the equation (\ref{deq}) on a pair of sets $A$ and $B$, which is clearly valid if one restricts the equation to their subsets. Hence, by Cauchy-Schwartz inequality and Lemma \ref{Glib} (i), we have
$$|\tilde A+\xi \tilde B| > \frac{q}{2}\left(\frac{|\tilde A||\tilde B|}{|A||B|}\right)^2 = \frac{q}{2}(uv)^2,$$
and the same estimate for all the four first sets on the list (\ref{sets}). Furthermore, if there is no doubling in the last five sets and all the nine sets are disjoint, the cardinality of the union of them all, using
 $$
 q<|A||B| = |\tilde A||\tilde B|+ |\tilde A||\bar B|+ |\bar A||\tilde B|+ |\bar A||\bar B|,
 $$ is greater than
$$
2q(uv)^2  + |A||B|( 2u(1-v) + 2v(1-u) + (1-u)(1-v))>q(2(uv)^2 + (u+v) - 3uv +1).
$$
It is easy to see that for $u,v\in (0,1)$,
$$
f(u,v) = 2(uv)^2 + (u+v) - 3uv\geq 0.
$$
Indeed, $f(u,v)$ is non-negative on the boundary of the above domain for $(u,v)$ and has a single critical point $u=v=\frac{1}{2}$ inside, where $f(\frac{1}{2},\frac{1}{2})=\frac{3}{8}>0.$ (At a critical point we have
$$
4uv^2 -3v+1 =0\qquad\mbox{and}\qquad u=v.
$$
The function $4v^3 -3v+1>0$ on $(0,1)$ attains its absolute minimum, equal to zero, at $v=\frac{1}{2}$.)

\medskip
Thus, assuming that all the sets (\ref{sets}) are disjoint and there is no doubling within the last five leads to an absurd statement that their union has size greater than $q$. The alternative to this is that there exists a non-trivial $x$ in some pair-wise intersection of the sets $C^\pm_\xi,-C^\pm_\xi$. Lemma \ref{pr2} now kicks in and completes the proof of Theorem \ref{th}. $\Box$

\end{document}